\documentclass{article}
\usepackage{amsmath}
\usepackage{amssymb}
\title{Lectures on Cyclotomic Hecke Algebras}
\author{Susumu Ariki}
\date{}
\begin{document}
\maketitle

\newtheorem{Theorem}{Theorem}[section]
\newtheorem{Corollary}[Theorem]{Corollary}
\newtheorem{Proposition}[Theorem]{Proposition}
\newtheorem{Definition}[Theorem]{Definition}
\newtheorem{Example}[Theorem]{Example}
\newtheorem{Lemma}[Theorem]{Lemma}
\newtheorem{Remark}[Theorem]{Remark}

\newcommand{\BM}{\left[\begin{array}}
\newcommand{\EM}{\end{array}\right]}

\section{Introduction}
\footnotetext{Hiroshima 7-9 July 1999, Durham 19-29 July 1999}
The purpose of these lectures is to introduce 
the audience to the theory of cyclotomic Hecke 
algebras of type $G(m,1,n)$. These algebras were introduced by 
the author and Koike, Brou\'e and Malle independently. 
As is well known, group rings of Weyl groups 
allow certain deformation. It is true for Coxeter groups, which 
are generalization of Weyl groups. These algebras are now known as 
(Iwahori) Hecke algebras. 

Less studied is its generalization to complex reflection 
groups. As I will explain later, this generalization is 
not artificial. The deformation of the group ring of the 
complex reflection group of type $G(m,1,n)$ is particularly 
successful. The theory uses many aspects of very modern 
development of mathematics: Lusztig and Ginzburg's geometric 
treatment of affine Hecke algebras, Lusztig's theory of 
canonical bases, Kashiwara's theory of global and crystal bases, 
and the theory of Fock spaces which arises from the study of 
solvable lattice models in Kyoto school. 

This language of Fock spaces is crucial in the theory 
of cyclotomic Hecke algebras. I would like to mention 
a little bit of history about Fock spaces in the context of 
representation theoretic study of solvable lattice models. 
For level one Fock spaces, it has origin in Hayashi's work. 
The Fock space we use is due to Misra and Miwa. 
For higher level Fock spaces, they appeared in work of 
Jimbo, Misra, Miwa and Okado, and Takemura and Uglov. 
We also note that Varagnolo and Vasserot's version of level one 
Fock spaces have straight generalization to 
higher levels and coincide with the Takemura and Uglov's one.
The Fock spaces we use are different from 
them. But they are essential in the proofs. 

Since the cyclotomic Hecke algebras contain the Hecke algebras 
of type A and type B as special cases, the theory of cyclotomic 
Hecke algebras is also useful to 
study the modular representation theory of finite classical 
groups of Lie type. 

I shall explain theory of Dipper and James, and 
its relation to our theory. The relevant Hecke algebras are 
Hecke algebras of type A. In this case, we have an alternative 
approach depending on the Lusztig's conjecture on quantum groups, 
by virtue of Du's refinement 
of Jimbo's Schur-Weyl reciprocity. 
Even for this rather well studied case, 
our viewpoint gives a new insight. 
This viewpoint first appeared in 
work of Lascoux, Leclerc and Thibon. 
This Fock space description 
looks quite different from the Kazhdan-Lusztig combinatorics, 
since it hides affine Kazhdan-Lusztig 
polynomials behind the scene. Inspired by this description, 
Goodman and Wenzl have 
found a faster algorithm to compute these polynomials. 
Leclerc and Thibon are key players 
in the study of this type A case. I also would like to mention 
Schiffman and Vasserot's work here, since it makes the relation of 
canonical bases between modified quantum algebras and quantized 
Schur algebras very clear. 

I will refer to work of Geck, Hiss, and Malle a little if time 
allows, since 
we can expect future development in this direction. It is 
relevant to Hecke algebras of type B. Finally, 
I will end the lectures with Brou\'e's famous dream. 

Detailed references can be found at the end of these lectures. 
The first three are for overview, and the rest are 
selected references for the lectures. [i-] implies a reference 
for the $i$ th lecture. 

\footnotetext{I would like to thank all the researchers involved in the 
development. Good interaction with German modular representation 
group (Geck, Hiss, Malle; Dipper), British combinatorial modular 
representation group (James, Mathas, Murphy), 
French combinatorics group (Lascoux, Leclerc, Thibon), 
modular representation group (Brou\'e, Rouquier; Vigneras), 
geometric representation group (Varagnolo, Vasserot, Schiffman) 
and Kyoto solvable lattice model group (Okado, Takemura, Uglov) 
has nourished the rapid development. We still have some problems 
to solve, and welcome young people who look for problems. 

I also thank Kashiwara, Lusztig, Ginzburg for their theories 
which we use. }


\section{Lecture One}
\subsection{Definitions}

Let ${\it k}$ be a field (or an integral domain in general). We define cyclotomic Hecke algebras of 
type $G(m,1,n)$ as follows. 

\begin{Definition}
Let $v_1,\dots,v_m, q$ be elements 
in ${\it k}$, and assume that $q$ is invertible. 
The Hecke algebra ${\cal H}_n(v_1,\dots,v_m;q)$ of type $G(m,1,n)$ 
is the ${\it k}$-algebra 
defined by the following relations for generators 
$a_i$ $(1\!\le\! i\!\le\! n)$. We often write ${\cal H}_n$ 
instead of ${\cal H}_n(v_1,\dots,v_m;q)$. If we want to 
make the base ring explicit, we write ${\cal H}_n/k$. 
\[
\begin{array}{c}
(a_1-v_1)\cdots(a_1-v_m)=0,\qquad (a_i-q)(a_i+1)=0\quad(i\ge 2)\\
\\
a_1a_2a_1a_2=a_2a_1a_2a_1, \qquad a_ia_j=a_ja_i\quad(j\ge i\!+\!2)\\
\\
a_ia_{i-1}a_i=a_{i-1}a_ia_{i-1}\;(3\!\le\! i\!\le\! n)
\end{array}
\]
The elements $L_i=q^{1-i}a_ia_{i-1}\cdots a_2a_1a_2\cdots a_i$ 
$(1\le i\le n)$ are called (Jucy-) Murphy elements or Hoefsmit 
elements. 
\end{Definition}

\begin{Remark}
Let $\hat H_n$ be the (extended) affine Hecke algebra associated 
with the general linear group over a non-archimedian field. 
For each choice of positive root system, we have Bernstein 
presentation of this algebra. Let $P={\mathbb Z}\epsilon_1+
\cdots+{\mathbb Z}\epsilon_n$ be the weight lattice as usual. 
We adopt "geometric choice" for the positive root system. 
Namely $\{\alpha_i:=\epsilon_{i+1}-\epsilon_i\}$ are simple roots. 
Let $S$ be the associated set of Coxeter generators (simple reflections). 
Then $\hat H_n$ has description via generators $X_\epsilon$ 
$(\epsilon\in P)$ and $T_s$ $(s\in S)$. We omit the description 
since it is well known. The following mapping gives rise a surjective 
algebra homomorphism from $\hat H_n$ to ${\cal H}_n$. 
\[
X_{\epsilon_i}\mapsto L_i, \qquad
T_{s_{\alpha_i}}\mapsto a_{i+1}
\]
This fact is the reason why we can apply Lusztig's theory to the 
study of cyclotomic Hecke algebras. Since the module theory 
for ${\cal H}_n$ has been developed by different methods, 
it has also enriched the theory of affine Hecke algebras. 
\end{Remark}

\begin{Remark}
Let $\zeta_m$ be a primitive $m$ th root of unity. If we 
specialize $q=1, v_i=\zeta_m^{i-1}$, we have the group ring 
of $G(m,1,n)$. $G(m,1,n)$ is the group of $n\times n$ permutation 
matrices whose non zero entries are allowed to be $m$ th roots of unity. 
Under this specialization, $L_i$ corresponds to the diagonal matrix 
whose $i$ th diagonal entry is $\zeta_m$ and whose remaining 
diagonal entries are $1$. We would like to stress two major 
differences between the group algebra and the deformed algebra 
${\cal H}_n$. 

\bigskip
\noindent
(1) $(L_i-v_1)\cdots(L_i-v_m)$ 
is not necessarily zero for $i>1$. 

\bigskip
\noindent
(2) If we consider 
the subalgebra generated by Murphy elements, its dimension is 
not $m^n$ in general. Further, the dimension depends on 
parameters $v_1,\dots,v_m,q$. 
\end{Remark}

Nevertheless, we have the following Lemma. $a_w$ is defined by 
$a_{i_1}\cdots a_{i_l}$ for a reduced word $s_{i_1}\cdots s_{i_l}$ 
of $w$. It is known that $a_w$ does not depend on the choice of 
the reduced word. 

\begin{Lemma}
$\{L_1^{e_1}\cdots L_n^{e_n}a_w|0\le e_i<m, w\in{\mathfrak S}_n\}$ 
form a basis of ${\cal H}_n$. 
\end{Lemma}
(How to prove) We consider ${\cal H}_n$ over an integral domain $R$, 
and show that $\sum RL_1^{e_1}\cdots L_n^{e_n}a_w$ is a two sided ideal. 
Then we have that these elements generate ${\cal H}_n$ as an 
$R$-module. To show that they are linearly independent, it is enough to 
take $R={\mathbb Z}[{\bf q},{\bf q}^{-1},{\bf v}_1,\dots,{\bf v}_m]$. 
In this generic parameter case, we embed the algebra into 
${\cal H}_n/
{\mathbb Q}({\bf q},{\bf v}_1,\dots,{\bf v}_m)$. Then we can 
construct enough simple modules to evaluate the dimension. $\blacksquare$

\bigskip
\noindent
An important property of ${\cal H}_n$ is the following. 

\begin{Theorem}[Malle-Mathas] Assume that $v_i$ are all invertible. 
Then ${\cal H}_n$ is a symmetric algebra.
\end{Theorem}
(How to prove) Since ${\cal H}_n$ is deformation of the group 
algebra of $G(m,1,n)$, we can define a length function 
$l(w)$ and $a_w$ for a reduced word of $w$. 
Unlike the Coxeter group case, $a_w$ does depend on the choice 
of the reduced word. Nevertheless, the trace function 
\[
tr(a_w)=\left\{\begin{array}{lr} 0 \quad(w\ne1)\\
                                 1 \quad(w=1)\end{array}\right.
\]
is well defined. $(u,v):=tr(uv)$ $(u,v\in{\cal H}_n)$ gives 
the bilinear form with the desired properties. $\blacksquare$

\begin{Remark}
We have defined deformation algebras for (not all but most of) other types of 
irreducible complex reflection groups by generators and 
relations. ($G(m,p,n)$: the author, 
other exceptional groups: Brou\'e and Malle.) 

The most natural definition of cyclotomic Hecke algebras is 
given by Brou\'e, Malle and Rouquir. It coincides with the 
previous definition in most cases. 
\end{Remark}

\bigskip
Let ${\cal A}$ be the hyperplane arrangement defined by complex 
reflections of $W$. For each ${\cal C}\in{\cal A}/W$, we can associate 
the order $e_{\cal C}$ of the cyclic group 
which fix a hyperplane in ${\cal C}$. Primitive idempotents of this 
cyclic group are denoted by $\epsilon_j(H)$ $(0\le j<e_{\cal C})$. 
We set ${\cal M}={\mathbb C}^n
\setminus\cup_{H\in{\cal A}}H$. 

\begin{Definition}
For each hyperplane $H$, 
let $\alpha_H$ be the linear form whose kernel is $H$. 
It is defined up to scalar multiple. We fix a set of 
complex numbers $t_{{\cal C},j}$. Then 
the following partial differential equation for 
${\mathbb C}W$-valued functions $F$ on ${\cal M}$ 
is called the (generalized) KZ equation. 
\[
\frac{\partial F}{\partial x_i}=\frac{1}{2\pi\sqrt{-1}}
\displaystyle\sum_{{\cal C}\in{\cal A}/W}
\sum_{j=0}^{e_{\cal C}-1}\sum_{H\in{\cal C}}
\frac{\partial(\log \alpha_H)}{\partial x_i}
t_{{\cal C},j}\epsilon_j(H)F
\]
\end{Definition}

\begin{Theorem}[Brou\'e-Malle-Rouquier] 
Assume that parameters are sufficiently generic. 
Let $B$ be the braid group attached to ${\cal A}$. Then the monodromy 
representation of $B$ with respect to the above KZ equation factors through 
a deformation ring of ${\mathbb C}W$. If $W=G(m,1,n)$ for 
example, it coincides with the cyclotomic Hecke algebra with 
specialized parameters. 
\end{Theorem}

\subsection{Representations}

If all modules are projective modules, we say that ${\cal H}_n$ 
is a semi-simple algebra, and call these representations {\bf ordinary 
representations}. We have 

\begin{Proposition}[Ariki(-Koike)] 
${\cal H}_n$ is semi-simple if and only if 
$q^iv_j-v_k$ $(|i|<n,j\ne k)$ and $1+q+\cdots+q^i$ $(1\le i<n)$ 
are all non zero. In this case, simple modules are parametrized 
by $m$-tuples of Young diagrams of total size $n$. For each 
$\lambda=(\lambda^{(m)},\dots,\lambda^{(1)})$, the 
corresponding simple module can be realized on the space 
whose basis elements are indexed by standard tableaux of shape $\lambda$. 
The basis elements are simultaneous eigenvectors of Murphy elements, 
and we have explicit matrix representation for generators $a_i$ 
$(1\le i\le n)$. 
\end{Proposition}

These represenations are called {\bf semi-normal form representations}. 

\noindent
Hence we have complete understanding of ordinary representations. 
If ${\cal H}_n$ is not semi-simple, representations are called 
{\bf modular representations}. A basic tool to get information for 
modular representations from ordinary ones is "reduction" 
procedure. 

\begin{Definition}
Let $(K,R,k)$ be a modular system. Namely, $R$ is a discrete 
valuation ring, $K$ is the field of fractions, and $k$ is the 
residue field. For an ${\cal H}_n/K$-module $V$, we take 
an ${\cal H}_n/R$-lattice $V_R$ and set $\overline V=V_R\otimes k$. 
It is known that $\overline V$ does depend on the choice of 
$V_R$, but the composition factors do not depend on the choice of 
$V_R$. The map between Grothendieck groups of finite dimentional 
modules given by 
\[
dec_{K,k}:K_0(mod\!\!-\!\!{\cal H}_n/K)\longrightarrow 
K_0(mod\!\!-\!\!{\cal H}_n/k)
\]
which sends $[V]$ to $[\overline V]$ is called a {\bf decomposition 
map}. Since Grothendieck groups have natural basis given by 
simple modules, we have the matrix representation of the 
decomposition map with respect to these bases. It is called 
the {\bf decomposition matrix}. The entries are called {\bf decomposition 
numbers}. 
\end{Definition}

In the second lecture, we also consider the decomposition map 
between Grothendieck groups of $KGL(n,q)\!\!-\!\!mod$ and 
$kGL(n,q)\!\!-\!\!mod$. 

\begin{Remark}
Decomposition maps are not necessarily surjective even after 
coefficients are extended to complex numbers. 
If we take $m=1,2$ and $q\in k$ 
to be zero, we have counter examples. These are called zero 
Hecke algebras, and studied by Carter. {\bf Note that we exclude 
the case $q=0$ in the definition}. In the case of group algebras, 
the theory of Brauer characters ensures that decomposition maps 
are surjective. 
\end{Remark}

In the case of cyclotomic Hecke algebras, we have the following 
result. 

\begin{Theorem}[Graham-Lehrer]
${\cal H}_n$ is a cellular algebra. In particular, the decomposition 
maps are surjective. 
\end{Theorem}

The notion of cellularity is introduced by Graham and Lehrer. 
It has some resemblance to the definition of quasi hereditary 
algebras. This is further pursued by K\"onig and Changchang Xi. 

In this lecture, we follow Dipper, James and Mathas' construction 
of Specht modules. We first fix notation. 

\bigskip
Let $\lambda=(\lambda^{(m)},\dots,\lambda^{(1)})$, 
$\mu=(\mu^{(m)},\dots,\mu^{(1)})$ be two $m$-tuples of 
Young diagrams. We say $\lambda$ dominates $\mu$ and 
write $\lambda\trianglerighteq\mu$ if 
\[
\displaystyle \sum_{j>k}|\lambda^{(j)}|+
\sum_{j=1}^l\lambda^{(k)}_j\ge
\sum_{j>k}|\mu^{(j)}|+
\sum_{j=1}^l\mu^{(k)}_j
\]
for all $k,l$. This partial order is called {\bf dominance order}. 

For each $\lambda=\left(\lambda^{(m)},\dots,\lambda^{(1)}\right)$, 
we set $a_k=n-|\lambda^{(1)}|-\cdots-|\lambda^{(k)}|$. 

\noindent
We have $n\ge a_1\ge\dots\ge a_l>0$ and $a_k=0$ 
for $k>l$ for some $l$. we denote $l$ by $l(a)$. 
For $a=(a_k)$, we denote by ${\mathfrak S}_a$ 
the set of permutations which preserve 
$\{1,\dots,a_l\},\dots,\{a_k+1,\dots,a_{k-1}\},\dots
\{a_1+1,\dots,n\}$. 
We also set
\[
\begin{array}{rl}
u_a=&(L_1-v_1)\cdots(L_{a_1}-v_1)\times(L_1-v_2)\cdots(L_{a_2}-v_2)
\times\cdots\\
&\cdots\times
(L_1-v_{l(a)})\cdots(L_{l(a)}-v_{l(a)})
\end{array}
\]

Let $t^\lambda$ be the canonical tableau. 
It is the standard tableau on which $1,\dots,n$ are filled 
in by the following rule; 

\noindent
$1,\dots,\lambda^{(m)}_1$ are written in the first row 
of $\lambda^{(m)}$; 
$\lambda^{(m)}_1+1,\dots,
\lambda^{(m)}_1+\lambda^{(m)}_2$ are written in the second row 
of $\lambda^{(m)}$; $\dots$; $|\lambda^{(m)}|+1,\dots,
|\lambda^{(m)}|+\lambda^{(m-1)}_1$ are written in the first row 
of $\lambda^{(m-1)}$; and so on. 

The row stabilizer of $t^\lambda$ 
is denoted by ${\mathfrak S}_\lambda$. We set
\[
x_\lambda=\displaystyle\sum_{w\in{\mathfrak S}_\lambda}a_w, 
\quad m_\lambda=x_\lambda u_a=u_ax_\lambda.
\]
Let $t$ be a standard tableau of shape $\lambda$. 
If the location of $i_k\in\{1,\dots,n\}$ in $t$ 
is the same as the location of $k$ in $t^\lambda$, 
We define $d(t)\in{\mathfrak S}_n$ by $k\mapsto i_k$ 
$(1\le k\le n)$. 

\begin{Definition}
Let $*:{\cal H}_n\rightarrow{\cal H}_n$ be the anti-involution 
induced by $a_i^*=a_i$. For each pair $(s,t)$ of standard 
tableaux of shape $\lambda$, we set $m_{st}=a_{d(s)}^*m_\lambda
a_{d(t)}$. 
\end{Definition}

\begin{Remark}
$\{m_{st}\}$ form a cellular basis of ${\cal H}_n$.
\end{Remark}

\begin{Proposition}[Dipper-James-Mathas]
Let $(K,R,k)$ be a modular system. We set 
${\cal I}_\lambda=\sum Rm_{st}$ where sum is over 
pairs of standard tableaux of shape strictly greater 
than $\lambda$ (with respect to the dominance order). 
Then ${\cal I}_\lambda$ is a two sided ideal of 
${\cal H}_n/R$.
\end{Proposition}
(How to prove) It is enough to consider straightening laws for 
elements $a_im_{st}$ and $m_{st}a_i$. We can then show that 
$m_{uv}$ appearing in the expression have greater shapes 
with respect to the dominance order. $\blacksquare$

\begin{Definition}
Set $z_\lambda=m_\lambda\;{\rm mod}\;{\cal I}_\lambda$. 
Then the submodule $S^\lambda=z_\lambda{\cal H}_n$ of 
${\cal H}_n/{\cal I}_\lambda$ is called a {\bf Specht module}. 
\end{Definition}

\begin{Theorem}

{\bf(Dipper-James-Mathas)}

\noindent
$\{z_\lambda a_{d(t)}|\;t:\mbox{standard of shape $\lambda$}\}$ 
form a basis of $S^\lambda$. 
\end{Theorem}
(How to prove) We can show 
by induction on the dominance order 
that these generate $S^\lambda$. Hence the collection 
of all these generate ${\cal H}_n$. Thus counting argument 
completes the proof. $\blacksquare$

\begin{Definition}
$S^\lambda$ is equipped with a bilinear form defined by
\[
\langle z_\lambda a_{d(t)},z_\lambda a_{d(s)}\rangle m_\lambda
=m_\lambda a_{d(s)}a_{d(t)}^*m_\lambda \;{\rm mod}\;{\cal I}_\lambda
\]
\end{Definition}

\begin{Theorem}

{\bf(General theory of Specht modules)}

\noindent
(1) $D^\lambda=S^\lambda/{\rm rad}\langle\hphantom{w},\hphantom{i}\rangle$ 
is absolutely irreducible or zero module. 
$\{D^\lambda\ne0\}$ form a complete set of simple ${\cal H}_n$-
modules.

\noindent
(2) Assume $D^\mu\ne0$ and $[S^\lambda:D^\mu]\ne0$. Then we have 
$\mu\trianglelefteq\lambda$. 
\end{Theorem}

\begin{Remark}
In the third lecture, we give a criterion for non vanishing 
of $D^\lambda$.
\end{Remark}

\begin{Theorem}[Dipper-Mathas]
Let $\{v_1,\dots,v_m\}=\sqcup_{i=1}^a S_i$ be the decomposition 
such that $v_j,v_k$ are in a same $S_i$ if and only if $v_j=v_kq^b$ 
for some $b\in{\mathbb Z}$. Then we have 
\[
mod\!\!-\!\!{\cal H}_n\simeq
\bigoplus_{n_1,\dots,n_a}mod\!\!-\!\!{\cal H}_{n_1}\boxtimes
\cdots\boxtimes mod\!\!-\!\!{\cal H}_{n_a}
\]
where ${\cal H}_n={\cal H}_n(v_1,\dots,v_m;q)$, 
${\cal H}_{n_i}={\cal H}_{n_i}(S_i;q)$, and the sum 
runs through $n_1+\cdots+n_a=n$. 
\end{Theorem}

Hence, it is enough to consider the case that $v_i$ are powers 
of $q$. 

\begin{Remark}
For the classification of simple modules, we can use 
arguments of Rogawski and Vigneras for the reduction 
to the case that $v_i$ are powers of $q$. 
Hence we do not need the above theorem 
for this purpose. 
\end{Remark}

\subsection{First application}
Let $k_q^\times=k^\times/\langle q\rangle$. We assume that $q\ne 1$, 
and denote the multiplicative order of $q$ by $r$. 
A {\bf segment} is a finite sequence of consecutive residue 
numbers which take values in ${\mathbb Z}/r{\mathbb Z}$. 
A {\bf multisegment} is a collection segments. Assume that 
a multisegment is given. Take a segment 
in the multisegment. 
By adding $i$ $(i\in{\mathbb Z}/r{\mathbb Z})$ 
to the entries of the segment simultaneously, 
we have a segment of 
shifted entries. 
If all of these $r$ segments appear in the given multisegment, 
we say that the given multisegment is {\bf periodic}. If it never 
happens for all segements in the multisegment, we say that 
the given multisegment is {\bf aperiodic}. We denote by 
${\cal M}_r^{ap}$ the set of aperiodic multisegments. 

\begin{Theorem}[Ariki-Mathas]
Simple modules over $\hat H_n/k$ are parametrized by 
\[
{\cal M}_r^{ap}(k)=\{\lambda:k_q^\times\rightarrow
{\cal M}_r^{ap}| \sum_{x\in k_q^\times}|\lambda(x)|=n\}
\]
\end{Theorem}
(How to prove) We consider a setting for reduction procedure, 
and show that a lower bound and an upper bound for the number 
of simple modules coincide. To achieve the lower bound, we use the integral 
module structure of the direct sum of Grothendieck groups 
of $proj\!\!-\!\!{\cal H}_n$ with respect to a Kac-Moody algebra 
action, which will be explained in the second lecture. 
The upper bound is achieved by cellularity. $\blacksquare$ 

\begin{Remark}
The lower bound can be achieved by a different method. 
This is due to Vigneras. 
\end{Remark}

Let $F$ be a nonarchimedian local field and assume that 
the residue field has characteristic different from the 
characteristic of $k$. We assume that $k$ is algebraically 
closed. We consider admissible $k$-representations 
of $GL(n,F)$. We take modular system $(K,R,k)$ and consider 
reduction procedure. 

\begin{Theorem}[Vigneras]
All cuspidal representations are obtained by reduction procedure. 
The admissible dual of $k$-representations is obtained from 
the classification of simple $\hat H_n/k$-modules. 
\end{Theorem}

Hence we have contribution to the last step of the classification. 

\begin{Remark}
Her method is induction from open compact groups and 
theory of minimal K-types. In the characteristic zero 
case, it is done by Bushnell and Kutzko. Considering 
$M:=ind_{G,K}(\sigma)$ where $(K,\sigma)$ is irreducible 
cuspidal distinguished K-type, she shows that 
${\rm End}_{kG}(M)$ is isomorphic to 
product of affine Hecke algebras, and $M$ 
satisfies the following hypothesis.
\end{Remark}

\medskip
"There exists a finitely generated projective module $P$ and 
a surjective homomorphism $\beta:P\rightarrow M$ such that 
${\rm Ker}(\beta)$ is ${\rm End}_{kG}(P)$-stable." 

\bigskip
Then the classification of simple $kG$-modules reduces to 
that of simple ${\rm End}_{kG}(M)$-modules. This simple 
fact is known as Dipper's lemma. 

\section{Lecture Two}
\subsection{Geometric theory}
Let ${\cal N}$ be the set of $n\times n$ nilpotent matrices, 
${\cal F}$ be the set of $n$-step complete flags in ${\mathbb C}^n$. 
We define the {\bf Steinberg variety} as follows. 
\[
Z=\{(N,F_1,F_2)\in{\cal N}\times
{\cal F}\times{\cal F}|\mbox{$F_1, F_2$ are $N$-stable}\}
\]
$G:=GL(n,{\mathbb C})\times{\mathbb C}^\times$ naturally acts 
on $Z$ via 
\[
(g,q)(N,F_1,F_2)=(q^{-1}Ad(g)N,gF_1,gF_2). 
\]
Let $K^G(Z)$ be the Grothendieck group of $G$-equivariant 
coherent sheaves on $Z$. It is an ${\mathbb Z}[{\bf q},{\bf q}^{-1}]$-
algebra via convolution product. 

\begin{Theorem}

{\bf(Ginzburg)}

\noindent
(1) We have an algebra isomorphism $K^G(Z)\simeq \hat H_n$. 

\noindent
(2) Let us consider a central character of the center 
${\mathbb Z}[X_{\epsilon_1}^\pm,
\dots,X_{\epsilon_n}^\pm]^{\mathfrak S_n}[{\bf q}^\pm]$ 
induced by $\hat s:X_{\epsilon_i}\mapsto 
\lambda_i$. By specializing the center 
via this linear character, we obtain a specialized affine 
Hecke algebra. Let $s$ be ${\rm diag}(\lambda_1,\dots,
\lambda_n)$. Then $H_*(Z^{(s,q)},{\mathbb C})$ equipped with 
convolution product is isomorphic to the specialized affine Hecke 
algebra. Here the homology groups are Borel-Moore homology groups, and 
$Z^{(s,q)}$ are fixed points of $(s,q)\in G$. 
\end{Theorem}

\begin{Remark}
All simple modules are obtained as simple modules of 
various specialized affine Hecke algebras. 
\end{Remark}

\begin{Theorem}

{\bf(Sheaf theoretic interpretation)}

Let $\tilde{\cal N}$ be 
$\{(N,F)\in{\cal N}\times{\cal F}|\mbox{$F$ is $N$-stable}\}$, 
$\mu:\tilde{\cal N}\rightarrow{\cal N}$ be the first projection. 
Then

\medskip
\noindent
(1) $H_*(Z^{(s,q)},{\mathbb C})\simeq 
Ext^*(\mu_*{\mathbb C}_{\tilde{\cal N}^{(s,q)}},
\mu_*{\mathbb C}_{\tilde{\cal N}^{(s,q)}})$. 

\noindent
(2) Let $\mu_*{\mathbb C}_{\tilde{\cal N}^{(s,q)}}=
\oplus_{\cal O}\oplus_{k\in{\mathbb Z}}L_{\cal O}(k)\otimes
IC({\cal O},{\mathbb C})[k]$. Then 
$L_{\cal O}:=\oplus_{k\in{\mathbb Z}} L_{\cal O}(k)$ is a simple 
$H_*(Z^{(s,q)},{\mathbb C})$-module or zero module. 
Further, non-zero ones form a complete set of simple 
$H_*(Z^{(s,q)},{\mathbb C})$-modules. If $q$ is not a root 
of unity, all $L_{\cal O}$ are non-zero. If $q$ is a primitive 
$r$ th root of unity, $L_{\cal O}\ne0$ if and only if 
$\cal O$ corresponds to a (tuple of) aperiodic multisegments taking residues 
in ${\mathbb Z}/r{\mathbb Z}$. 
\end{Theorem}

In the above theorem, the orbits run through orbits consisting 
of isomorphic representations of a quiver, 
which is disjoint union of infinite line quivers 
or cyclic quivers of length $r$. 
The reason 
is that ${\cal N}^{(s,q)}$ is the set of nilpotent matrices 
$N$ satisfying $sNs^{-1}=qN$, which can be identified with 
representations of a quiver via considering eigenspaces of $s$ 
as vector spaces on nodes and $N$ as linear maps on arrows. 
This is the key fact which relates the affine quantum algebra 
of type $A_\infty$, $A_{r-1}^{(1)}$ and representations of cyclotomic 
Hecke algebras. 

\begin{Definition}
Let ${\cal C}_n$ be the full subcategory of $mod\!\!-\!\!\hat H_n$ 
whose objects are modules which have central character $\hat s$ 
with all eigenvalues of $s$ being powers of $q$. Set 
$c_n=X_{\epsilon_1}+\cdots+X_{\epsilon_n}$. We denote by 
$P_{c_n,\lambda}(-)$ the exact functor taking generalized 
eigenspaces of eigenvalue $\lambda$ with respect to $c_n$. 
We then set
\[
i-Res(M)=\bigoplus_{\lambda\in k}
P_{c_{n-1},\lambda-q^i}\left(Res_{\hat H_{n-1}}^{\hat H_n}\left(
P_{c_n,\lambda}(M)\right)\right)
\]
This is an exact functor from ${\cal C}_n$ to ${\cal C}_{n-1}$. 
We set $U_n={\rm Hom}_{\mathbb C}(K_0({\cal C}_n),{\mathbb C})$, 
$f_i=(i-Res)^T:U_{n-1}\rightarrow U_n$. 
\end{Definition}

I shall give some historical comments here. 
The motivation to introduce these definitions was 
Lascoux-Leclerc-Thibon's 
observation that Kashiwara's global basis on level one modules 
computes the decomposition numbers of Hecke algebras of type A 
over the field of complex numbers. The above notions for 
affine Hecke algebras and cyclotomic Hecke algebras were first 
introduced by the author in his interpretation of 
Fock spaces and action of 
Chevalley generators in LLT observation 
into (graded dual of) Grothendieck groups of 
these Hecke algebras and 
$i$-restriction and $i$-induction operations. 
This is the starting point of a new point of view 
on the representation theory of affine Hecke algebras 
and cyclotomic Hecke algebras. As I will explain below, 
it allows us to give a new application of Lusztig's 
canonical basis. It triggered intensive studies of 
canonical bases on Fock spaces. These are carried out 
mostly in Paris and Kyoto. On the other hand, 
the research on cyclotomic Hecke algebras are mostly 
lead by Dipper, James, Mathas, Malle and the author. 
In the third lecture, these two will be combined to prove 
theorems on Specht module theory of cyclotomic Hecke 
algebras. 

We now state a key proposition necessary for the proof 
of the next theorem. In the top row of the diagram, we allow 
certain infinite sum in $U({\mathfrak g}(A_\infty))$ 
in accordance with infinite sum in $U_n$. Note that 
we do not have infinite sum in the bottom row. 

\begin{Proposition}[Ariki]
There exists a commutative diagram 
\[
\begin{array}{ccc}
U^-({\mathfrak g}(A_\infty))&\simeq&
\bigoplus_{n\ge0}U_n/{\bf q}\\
\uparrow& & \uparrow\\
U^-({\mathfrak g}(A_{r-1}^{(1)}))&\simeq&
\bigoplus_{n\ge0}U_n/q\!\!=\!\!\sqrt[r]{1}
\end{array}
\]
such that the left vertical arrow is inclusion, the 
right vertical arrow is induced by specialization 
${\bf q}\rightarrow q$, and the bottom horizontal 
arrow is an $U^-({\mathfrak g}(A_{r-1}^{(1)}))$-module 
isomorphism. Under this isomorphism, 
canonical basis elements of $U^-({\mathfrak g}(A_{r-1}^{(1)}))$ 
map to dual basis elements of $\{[\mbox{simple module}]\}$. 
\end{Proposition}
(How to prove) We firstly construct the upper horizontal arrow 
by using PBW-type basis and dual basis of 
$\{[\mbox{standard module}]\}$ of 
affine Hecke algebras. Here we use Kazhdan-Lusztig induction 
theorem. We also use restriction rule for Specht modules. 
We then appeal to folding argument. 
On the left hand side, we consider this folding in geometric 
terms. Since only short explanation was supplied in my original 
paper, I refer to Varagnolo-Vasserot's argument instead for 
this part. We then use 

\medskip
[standard module:simple module]=[canonical basis:PBW-type basis]

\medskip
\noindent
which is a consequence of the Ginzburg's theorem stated above. $\blacksquare$

We now turn to the cyclotomic case. In this case, we can 
consider not only negative part of Kac-Moody algebra, but 
the action of the whole Kac-Moody algebra. 

\begin{Definition}
Assume that $v_i=q^{\gamma_i}$ $(1\le i\le m)$ and $q=\sqrt[r]{1}$. 
We set 
\[
V_n={\rm Hom}_{\mathbb C}(K_0(mod\!\!-\!\!{\cal H}_n),
{\mathbb C}), \quad V=\oplus_{n\ge0}V_n. 
\]
We define $c_n=L_1+\cdots+L_n$. Then we can define 
\[
\begin{array}{c}
i-Res(M)=\displaystyle\bigoplus_{\lambda\in k}
P_{c_{n-1},\lambda-q^i}\left(Res_{{\cal H}_{n-1}}^{{\cal H}_n}\left(
P_{c_n,\lambda}(M)\right)\right),
\\
i-Ind(M)=\displaystyle\bigoplus_{\lambda\in k}
P_{c_{n+1},\lambda+q^i}\left(Ind_{{\cal H}_n}^{{\cal H}_{n+1}}\left(
P_{c_n,\lambda}(M)\right)\right).
\end{array}
\]
These are exact functors and we can define
\[
\begin{array}{rl}
e_i=(i-Ind)^T&: V_{n+1}\rightarrow V_n \\
f_i=(i-Res)^T&: V_{n-1}\rightarrow V_n
\end{array}
\]
\end{Definition}

\begin{Definition}
Let ${\cal F}=\oplus{\mathbb C}\lambda$ be a based vector space 
whose basis elements are $m$-tuples of Young diagrams 
$\lambda=(\lambda^{(m)},\dots,\lambda^{(1)})$. 

Assume that $\gamma_i\in{\mathbb Z}/r{\mathbb Z}$ $(1\le i\le m)$ 
are given. We introduce 
the notion of residues of cells as follows: 
Take a cell in $\lambda$. If the cell is located on 
the $(i,j)$th entry of $\lambda^{(k)}$, we say that the cell 
has residue $-i+j+\gamma_k\in{\mathbb Z}/r{\mathbb Z}$. 
Once residues are defined, we can speak of removable $i$-nodes 
and addable $i$-nodes on $\lambda$: Convex corners of $\lambda$ 
with residue $i$ are called {\bf removable $i$-nodes}. Concave 
corners of $\lambda$ with residue $i$ are called 
{\bf addable $i$-nodes}. 

We define operators $e_i$ and $f_i$ by 
$e_i\lambda$ (resp.$f_i\lambda)$ being the sum of all 
$\mu$'s obtained from $\lambda$ by removing (resp.adding) 
a removable (resp.addable) $i$-node. We can extend this 
action to make ${\cal F}$ an integrable ${\mathfrak g}(A_{r-1}^{(1)})$-
module. (If $r=\infty$, we consider $A_{r-1}^{(1)}$ as $A_\infty$.) 

We call ${\cal F}$ the {\bf combinatorial Fock space}. Note that 
the action of the Kac-Moody algebra 
depends on $(\gamma_1,\dots,\gamma_m;r)$. 
\end{Definition}

\begin{Theorem}[Ariki]
\label{LLT theorem} We assume $v_i=q^{\gamma_i}$ $(1\le i\le m)$, 
$q=\sqrt[r]{1}\in{\mathbb C}$. We set $\Lambda=\sum_{i=1}^m
\Lambda_{\gamma_i}$. Then we have the following. 

\bigskip
\noindent
(1) $L(\Lambda)\simeq V=U({\mathfrak g}(A_{r-1}^{(1)}))\emptyset
\subset{\cal F}$. 

\noindent
(2) Through this isomorphism, canonical basis elements of 
$L(\Lambda)$ are identified with dual basis elements of 
simple modules, and the embedding to ${\cal F}$ 
is identified with the transpose of the decomposition map. 
\end{Theorem}
(How to prove) We first consider reduction procedure from 
semi-simple ${\cal H}_n/K$ to ${\cal H}_n/k$. Note that 
this is not achieved by ${\bf v}_i={\bf q}^{\gamma_i}$ 
and ${\bf q}$ to $q$. Then $V/K$ can be identified with 
${\cal F}$. We then consider 
\[
\begin{array}{cccc}
U^-({\mathfrak g}(A_{r-1}^{(1)}))\emptyset&\simeq&
V & \subset {\cal F}\\
\uparrow& & \uparrow &\\
U^-({\mathfrak g}(A_{r-1}^{(1)}))&\simeq&
\oplus\;U_n/q\!\!=\!\!\sqrt[r]{1}& 
\end{array}
\]
Then the previous proposition and integrality of $\cal F$ 
prove the theorem. $\blacksquare$

\begin{Remark}
The theorem says that we have a new application of 
Lusztig's canonical bases, which is similar to the application 
of Kazhdan-Lusztig bases of Hecke algebras to Lie algebras 
(Kazhdan-Lusztig conjecture) and quantum algebras (Lusztig conjecture). 
It is interesting to observe that the roles of quantum algebras and 
Hecke algebras are interchanged: in Lusztig's conjecture, 
Kazhdan-Lusztig bases of Hecke algebras describe 
decomposition numbers of quantum algebras at roots of unity; 
in our case, canonical bases of quantum affine algebras on integrable modules 
describe decomposition numbers of cyclotomic Hecke algebras at 
roots of unity. Previously, 
a positivity result was the only application of canonical bases. 

The fact that affine Kazhdan-Lusztig 
polynomials appear in geometric construction 
of quantum algebras and affine Hecke algebras was known to 
specialists. What was new for affine Hecke algebras 
is the above proposition, particularly its formulation in terms 
of Grothendieck groups of affine Hecke algebras. 

For canonical bases on integrable modules, the theorem was 
entirely new, since no one knew 
the "correct" way of 
taking quotients of affine Hecke algebras to get the 
similar Grothendieck group description of canonical bases on 
integrable modules. It was just after 
cyclotomic Hecke algebras were introduced. 
\end{Remark}

\begin{Remark}
Let $(K,R,k)$ be a modular system. If we take semi-perfect $R$, 
we can identify $V$ with $\oplus_{n\ge0}K_0(proj\!\!-\!\!
{\cal H}_n)$, the transpose of the decomposition map with 
the map induced by lifting idempotents, the dual basis 
elements of simple modules with principal indecomposable 
modules, respectively. 
Here $proj\!\!-\!\!{\cal H}_n$ denotes the category 
of finite dimentional projective ${\cal H}_n$-modules. 
We often use this description since it is more appealing. 
\end{Remark}

\begin{Remark}
If $m=1$, namely the Hecke algebra has type A, we have another way 
to compute decomposition numbers. Let us consider Jimbo's Schur-
Weyl reciprocity. It has refinement by Du, and can be 
considered with specialized parameters. 
Let us denote the dimension of the natural representation 
by $d$, the endomorphism ring 
${\rm End}_{{\cal H}_n}(V^{\otimes n})$ by ${\cal S}_{d,n}$. 
This endomorphism ring is called 
{\bf the $q$-Schur algebra}. Note that Schur functors embed 
the decomposition numbers of Hecke algebras into those of 
$q$-Schur algebras. 
Then Du's result implies that 
the decomposition numbers of Hecke algebras 
are derived from those of 
quantum algebras $U_q({\mathfrak gl}_d)$ with 
$q=\sqrt[r]{1}$. 
There is a closed formula 
for decomposition numbers {\rm[Weyl module:simple module]} of 
quantum algebras at a root of unity: these are 
values at $1$ of parabolic Kazhdan-Lusztig polynomials for (extended) affine 
Weyl groups of type $A$. This formula is known as 
the Lusztig conjecture for quantum algebras. 
(This is a theorem of Kazhdan-Lusztig+Kashiwara-Tanisaki. 
There is another approach for this $m=1$ case. This is 
due to Varagnolo-Vasserot and Schiffman.) 
\end{Remark}

\begin{Remark}
The introduction of combinatorial Fock spaces is due to 
Misra, Miwa and Hayashi, as I stated in the introduction. 
We will return to 
their work on ${\bf v}$-deformed Fock spaces in the third lecture. 
\end{Remark}

\subsection{Algorithms}
For the case $m=1$, we have four algorithms to compute 
decomposition numbers. These are LLT algorithm, LT algorithm, 
Soergel algorithm, and modified LLT algorithm. For general $m$, 
we have Uglov algorithm. 

\medskip
\noindent
(1) LLT algorithm 

This is due to Lascoux, Leclerc and Thibon. It is based on 
theorem \ref{LLT theorem}. Basic idea is to construct 
"ladder decompostion" of restricted Young diagrams. Then 
it produces basis $\{A(\lambda)\}$ of the level one module 
$L(\Lambda_0)$. (I will show an example in the lecture. 
This is a very simple procedure.)

Once $\{A(\lambda)\}$ is given, we can determine canonical basis elements 
$G(\lambda)$ recursively. We set 
\[
G(\lambda)=A(\lambda)-\sum_{\mu\triangleright\lambda}
c_{\lambda\mu}(v)G(\mu), 
\]
and find $c_{\lambda\mu}(v)$ by the following condition.
\[
c_{\lambda\mu}(v^{-1})=c_{\lambda\mu}(v), \quad 
G(\lambda)\in \lambda+\sum_{\mu\triangleright\lambda}
v{\mathbb Z}[v]\mu
\]
Note that we follow the convention that restricted partitions 
form a basis of $L(\Lambda_0)$. 

\begin{Remark}
By a theorem of Leclerc, we can also compute decomposition 
numbers of $q$-Schur algebras by using those of Hecke algebras. 
\end{Remark}

\medskip
\noindent
(2) LT algorithm

This is based on Leclerc-Thibon's involution and Varagnolo-
Vasserot's reformulation of Lusztig conjecture. It has an 
advantage that we directly compute all decompoition 
numbers of $q$-Schur algebras. 

We use fermionic description of the Fock space. Then 
a simple procedure on basis elements and straightening laws 
define bar operation 
on the Fock space. We then compute canonical basis elements 
by the characterization 
\[
\overline{G(\lambda)}=G(\lambda), \quad
G(\lambda)\in \lambda+\sum_{\mu\triangleright\lambda}
v{\mathbb Z}[v]\mu
\]

\medskip
\noindent
(3) Soergel algorithm

It is reformulation of Kazhdan-Lusztig 
algorithm for parabolic Kazhdan-Lusztig polynomials. 
Let ${\cal A}^+$ be the set of alcoves 
in the positive Weyl chamber. We consider vector 
space with basis $\{(A)\}_{A\in{\cal A}^+}$. For each 
simple reflection $s$, we denote by $As$ the adjacent 
alcove obtained by the reflection. 
The Bruhat order determines partial order on ${\cal A}^+$. 
Let $C_s$ be the Kazhdan-Lusztig element corresponding to 
$s$ (we use $(T_s-v)(T_s+v^{-1})=0$ as a defining relation here). 
Then the action of $C_s$ on this space is given by 
\[
(A)C_s=\left\{\begin{array}{ll}
(As)+v(A)\quad&(As\in{\cal A}^+,As>A)\\
(As)+v^{-1}(A)&(As\in{\cal A}^+,As<A)\\
0 &(else)\end{array}\right.
\]

We determine Kazhdan-Lusztig basis elements $G(A)$ recursively. 
For $A\in{\cal A}^+$, we take $s$ such that $As<A$. Then
we find 
\[
G(A)=G(As)C_s-\sum_{B<A}c_{A,B}(v)G(B)
\]
by the condition
\[
c_{A,B}(v^{-1})=c_{A,B}(v),\quad 
G(A)\in (A)+\sum_{B<A}v{\mathbb Z}[v](B)
\]

\medskip
\noindent
(4) modified LLT algorithm

This is an algorithm which improves LLT algorithm. The idea 
is not to start from the empty Young diagram. This is due to 
Goodman and Wenzl. Their experiment shows that Soergel's is 
better than LLT, and modified LLT is much faster than both. 

\medskip
\noindent
(5) Uglov algorithm

This is generalization of LT algorithm, and it uses 
the higher level Fock space introduced by Takemura and Uglov. 

\subsection{Second application}

Let us return to the $q$-Schur algebra. We summarize 
the previous explanation as follows. 

\begin{Theorem}
If $q\ne 1$ is a root of unity in a field of characteristic zero, 
the decomposition numbers of the $q$-Schur algebra are 
computable. 
\end{Theorem}

\begin{Corollary}[Geck]
Let $k$ be a field. We consider the $q$-Schur algebra over 
$k$. If the characteristic of $k$ is 
sufficiently large, the decomposition numbers of 
the $q$-Schur algebra over $k$ are computable. 
Note that we do not exclude $q=1$ here. 
\end{Corollary}

It has application to the modular representation 
theory of $GL(n,q)$. Let 
$q$ be a power of a prime $p$, the characteristic 
of $k$ be $l\ne p$. We assume that $k$ is algebraically 
closed. This case is called non-describing 
characteristic case. We want to study $K_0(kGL(n,q)\!\!-\!\!mod)$. 

\begin{Theorem}[Dipper-James]
Assume that the 
decomposition numbers of $q^a$-Schur algebras over $k$ 
for various $a\in{\mathbb Z}$ are known. Then 
the decomposition numbers of $GL(n,q)$ in 
non-describing characteristic case are computable. 
\end{Theorem}

We explain how to compute the decomposition numbers 
of $G:=GL(n,q)$. Let $(K,R,k)$ be an $l$-modular system. 
James has constructed Specht modules for $RG$. We denote 
them by $\{S_R(s,\lambda)\}$. $s$ is a semi-simple element of 
$G$. If the degree of $s$ over ${\mathbb F}_q$ is $d$,  
$\lambda$ run through partitions of size $n/d$. 

\medskip
\noindent
(1) A complete set of simple $KG$-modules is given by 
\[
\left\{ R^G\left(
\underset{1\le i\le N}{\boxtimes} S_K(s_i,\lambda^{(i)})\right)\bigm|
\sum d_i|\lambda^{(i)}|=n \right\}
\]
where $R^G(-)$ stands for Harish-Chandra induction, $d_i$ 
is the degree of $s_i$, and 
$\{s_1,\dots,s_N\}$ run through sets of distinct 
semi-simple elements. We use 
Dipper-James' formula 
\[
[S_k(s,\lambda):D_k(s,\mu)]=d_{\lambda'\mu'}
\]
where $d_{\lambda'\mu'}$ is a decomposition number of the 
$q^d$-Schur algebra. Then we rewrite 
$R^G\left(\underset{1\le i\le N}{\boxtimes} S_k(s_i,\lambda^{(i)})\right)$ 
into sum of 
$R^G\left(\underset{1\le i\le N}{\boxtimes} D_k(s_i,\mu^{(i)})\right)$. 

\medskip
\noindent
(2) Let $t_i$ be the $l$-regular part 
of $s_i$, $a_i$ be the degree of $t_i$, $\nu^{(i)}$ 
be the Young diagram obtained from $\mu^{(i)}$ 
by multiplying all columns 
by $d_i/a_i$. Then we have $D_k(s_i,\mu^{(i)})\simeq
D_k(t_i,\nu^{(i)})$. This is also due to Dipper and James. 
Thus we can rewrite 
$R^G\left(\underset{1\le i\le N}{\boxtimes} D_k(s_i,\mu^{(i)})\right)$ 
into 
$R^G\left(\underset{1\le i\le N}{\boxtimes} D_k(t_i,\nu^{(i)})\right)$. 
Assume that $t_i=t_j$. Then we use the inverse of 
the decomposition matrices of $q^{a_i}$-Schur algebras of 
rank $d_ik_i$ and $d_jk_j$ to 
describe $D_k(t_i,\nu^{(i)})\boxtimes D_k(t_j,\nu^{(j)})$ as 
an alternating sum of 
$S_k(t_i,\eta^{(i)})\boxtimes S_k(t_j,\eta^{(j)})$. Then 
the Harish-Chandra induction of this module is explicitly 
computable by using Littlewood-Richardson rule. We use the 
decomposition matrix of the $q^{a_i}$-Schur algebra of rank 
$d_ik_i+d_jk_j$ to rewrite it again into the sum of 
$R^G\left(\underset{1\le i\le N'}{\boxtimes} D_k(t_i,\kappa^{(i)})\right)$. 
Continuing this procedure, we reach the case that 
all $t_i$ are mutually distinct. 

\medskip
\noindent
(3) The final result of the previous step already gives the answer 
since the following set is a complete set of simple $kG$-modules. 
\[
\left\{ R^G\left(\underset{1\le i\le N'}{\boxtimes} D_k(t_i,\kappa^{(i)})\right)\bigm|
\sum a_i|\kappa^{(i)}|=n \right\}
\]
where $\{t_1,\dots,t_{N'}\}$ run through sets of distinct 
$l$-regular semi-simple elements.

\section{Lecture Three}
\subsection{Specht modules and ${\bf v}$-deformed Fock spaces}
We now ${\bf v}$-deform the setting we have explained in the second 
lecture. The view point which has emerged is that behind 
the representation theory 
of cyclotomic Hecke algebras, there is the same crystal structure as 
integrable modules over quantum algebras of type 
$A_{r-1}^{(1)}$, and this crystal structure is induced by 
canonical bases of integrable modules. As a corollary to 
this viewpoint, Mathas and the author have parametrized 
simple ${\cal H}_n$-modules over an arbitrary field using 
crystal graphs. 
Since the canonical basis is defined in the ${\bf v}$-deformed setting, 
It further lead to intensive study of canonical bases on 
various ${\bf v}$-deformed Fock spaces. 

The purpose of the third lecture is to show the compatibility 
of this crystal structure with Specht module theory. The above 
mentioned studies on canonical bases on ${\bf v}$-deformed 
Fock spaces are essential in the proof. 

Before going to this main topic, I shall mention related work 
recently done in Vazirani's thesis. This can be understood in the 
above context. As I have explained in the second lecture, this viewpoint 
has origin in Lascoux, Leclerc and Thibon's work, which I would like 
to stress here again. 

\begin{Theorem}[Vazirani-Grojnowski]
Let $\tilde{e}_i(M)=soc(i-Res(M))$. If $M$ is irreducible, 
then $\tilde{e}_i(M)$ is irreducible or zero module. 
\end{Theorem}

The case $m=1$ is included in Kleshchev and Brundan's modular 
branching rule. It is natural to think that the socle series 
would explain the canonical basis in the crystal structure. 
This observation was first noticed by Rouquier as was explained 
in \cite{LLT}, and adopted in this Vazirani's thesis. 

We now start to explain how Specht module theory fits in 
the description of higher level Fock spaces. 

Let ${\mathcal F}_{\bf v}=\oplus {\mathbb C}({\bf v})\lambda$ be 
the ${\bf v}$-deformed Fock space. It has $U_{\bf v}({\mathfrak g}
(A_{r-1}^{(1)}))$-module structure which is deformation 
of $U({\mathfrak g}(A_{r-1}^{(1)}))$-module structure on 
${\cal F}$. To explain it, we introduce notation. 

Let $x$ be a cell on $\lambda=(\lambda^{(m)},\dots,\lambda^{(1)})$. 
Assume that it is the $(a,b)$th cell of $\lambda^{(c)}$. 
We say that a cell is {\bf above} $x$ if it is on $\lambda^{(k)}$ 
for some $k>c$, or if it is on $\lambda^{(c)}$ and the row number is 
strictly smaller than $a$. We denote the set of 
addable (resp. removable) $i$-nodes of $\lambda$ 
which are above 
$x$ by $A_i^a(x)$ (resp. $R_i^a(x)$). In a similar way, we say that 
a cell is {\bf below} $x$ if it is on $\lambda^{(k)}$ 
for some $k<c$, or if it is on $\lambda^{(c)}$ and the row number is 
strictly greater than $a$. We denote 
the set of addable (resp. removable) $i$-nodes of $\lambda$ 
which are 
below $x$ by $A_i^b(x)$ (resp. $R_i^b(x)$). 
The set of all addable (resp. removable) $i$-nodes of 
$\lambda$ is 
denoted by $A_i(\lambda)$ (resp. $R_i(\lambda)$). 
We then set 
\[
\begin{array}{c}
N_i^a(x)=|A_i^a(x)|-|R_i^a(x)|,\quad
N_i^b(x)=|A_i^b(x)|-|R_i^b(x)|\\
\\
N_i(\lambda)=|A_i(\lambda)|-
|R_i(\lambda)|
\end{array}
\]
We denote the number of all $0$-nodes in 
$\lambda$ by $N_d(\lambda)$. Then 
the $U_{\bf v}({\mathfrak g}(A_{r-1}^{(1)}))$-module structure given to 
${\mathcal F}_{\bf v}$ is as follows. 
\[
\begin{array}{c}
e_i{\lambda}=\displaystyle\sum_{\lambda/\mu=\fbox{i}}
{\bf v}^{-N_i^a(\lambda/\mu)}\mu, \quad
f_i{\lambda}=\displaystyle\sum_{\mu/\lambda=\fbox{i}}
{\bf v}^{N_i^b(\mu/\lambda)}\mu \\
\\
{\bf v}^{h_i}\lambda={\bf v}^{N_i(\lambda)}\lambda, \quad
{\bf v}^d\lambda={\bf v}^{-N_d(\lambda)}\lambda
\end{array}
\]
This action is essentially due to Hayashi. 

Set $V_{\bf v}=U_{\bf v}({\mathfrak g}(A_{r-1}^{(1)}))\emptyset$. It is 
considered as the ${\bf v}$-deformed space of 
$V=\oplus_{n\ge0}K_0(proj\!\!-\!\!{\cal H}_n)$. 

\begin{Remark}
If we apply a linear map $(\lambda^{(m)},\dots,\lambda^{(1)}) 
\mapsto ({\lambda^{(1)}}',\dots,{\lambda^{(m)}}')$, we have 
Kashiwara's lower crystal base which is compatible 
with his coproduct $\Delta_-$. 

On the other hand, if an anti-linear map 
$(\lambda^{(m)},\dots,\lambda^{(1)})
\mapsto ({\lambda^{(m)}}',\dots,{\lambda^{(1)}}')$ is applied, 
we have Lusztig's basis at $\infty$ which is compatible 
with his coproduct. We denote it by 
${\cal F}_{{\bf v}^{-1}}^{-\gamma}$. 
\end{Remark}

Set $L=\oplus {\mathbb Q}[{\bf v}]_{({\bf v})}\lambda$ 
and $B=\{\lambda\;\mbox{mod}\;{\bf v}\}$. Then it is known that 
$(L,B)$ is a crystal base of ${\cal F}_{\bf v}$. We nextly set 
$L_0=V_{\bf v}\cap L$, and $B_0=(L_0/{\bf v}L_0)\cap B$. Then 
general theory concludes that $(L_0,B_0)$ is a crystal base of 
$V_{\bf v}$. 

\begin{Definition}
We say that $\lambda$ is $(\gamma_1,\dots,\gamma_m;r)$-
Kleshchev if 
$\lambda\;\mbox{mod}\;{\bf v}\in B_0$. We often drop 
parameters and simply says $\lambda$ is Kleshchev. 
\end{Definition}

It has the following combinatorial definition. We say that 
a node on $\lambda$ is {\bf good} if there is 
$i\in{\mathbb Z}/r{\mathbb Z}$ such that if we read 
addable $i$-nodes (write A in short) and removable $i$-nodes 
(write R in short) from the 
top row of $\lambda^{(m)}$ to the bottom row of 
$\lambda^{(1)}$ and do RA deletion as many as possible, 
then the node sits in the left end of the remaining R's. 
(I will give an example in the lecture.) 

\begin{Definition}
$\lambda$ is called $(\gamma_1,\dots,\gamma_m;r)$-
{\bf Kleshchev} if there is a standard 
tableau $T$ of shape $\lambda$ such that for all $k$, 
$\fbox{k}$ is a good 
node of the subtableau $T_{\le k}$ which consists of 
nodes $\fbox{1},\dots,\fbox{k}$ by definition. 
\end{Definition}

\begin{Theorem}[Ariki]
We assume that $v_i=q^{\gamma_i},q=\sqrt[r]{1}$. Then 
$D^\lambda\ne0$ if and only if $\lambda$ is 
$(\gamma_1,\dots,\gamma_m;r)$-Kleshchev. 
\end{Theorem}
(How to prove) We show that canonical basis elements 
$G(\lambda)$ ($\lambda$=Kleshchev) have the form
\[
G(\lambda)=\lambda+\sum_{\mu\triangleright\lambda}c_{\lambda\mu}({\bf v})\mu
\]
On the other hand, the Specht module theory constructed by Dipper-James-Mathas 
shows that the principal indecomposable module 
$P^\lambda$ for $D^\lambda\ne0$ has the form
\[
P^\lambda=S^\lambda+\sum_{\mu\triangleright\lambda}m_{\lambda\mu}S^\mu
\]
Comparing these, and recalling that $\lambda\in{\cal F}$ is identified 
with $S^\lambda$, we have the result. $\blacksquare$

\bigskip
To know the form of $G(\lambda)$, we have to understand 
higher level ${\bf v}$-deformed Fock spaces. 

\begin{Definition}
Take $\gamma=(\gamma_1,\dots,\gamma_m)
\in({\mathbb Z}/r{\mathbb Z})^m$. If 
$\tilde\gamma=(\tilde\gamma_1,\dots,\tilde\gamma_m)\in{\mathbb Z}^m$ 
satisfies $\tilde\gamma_k\;\mbox{mod}\;r=\gamma_k$ 
for all $k$, we say that $\tilde\gamma$ is a {\bf lift} of 
$\gamma$. 
\end{Definition}

\begin{Theorem}[Takemura-Uglov]
For each $\tilde\gamma\in{\mathbb Z}^m$, we can construct higher level 
${\bf v}$-deformed Fock space, whose underlying space is 
the same as ${\cal F}_{\bf v}$. 
\end{Theorem}

It has geometric realization due to Varagnolo and Vasserot. 
For reader's convenience, I also add it here. 
Let $V$ be a ${\mathbb Z}$-graded ${\mathbb C}$-vector space 
whose dimension type 
is $(d_i)_{i\in{\mathbb Z}}$. We denote by $\overline{V}$ 
the ${\mathbb Z}/r{\mathbb Z}$-graded vector space defined by 
$\overline V_{\bar i}=\oplus_{j\in\bar i}V_j$. 
We set $\overline{V}_{j\ge i}=\oplus_{j\ge i}V_j$. Let 
\[
E_V=\underset{i\in{\mathbb Z}}{\oplus}{\rm Hom}_{\mathbb C}
(V_i,V_{i+1}), \quad 
E_{\overline V}=\underset{\bar i\in{\mathbb Z}/
r{\mathbb Z}}{\oplus}{\rm Hom}_{\mathbb C}
(V_{\bar i},V_{\overline{i+1}}). 
\]
and define $E_{\overline V,V}=\{x\in \overline V|\;x(\overline{V}_{j\ge i})
\subset\overline{V}_{j\ge i}\}$. Then we have a natural diagram 
\[
E_V \stackrel{\kappa}{\leftarrow} 
E_{\overline V,V} \stackrel{\iota}{\rightarrow} E_V
\]
We consider $\gamma_d:=\kappa_!\iota^*[{\rm shift}]$. Then it defines 
a map from the derived category which is used to construct 
$U_{\bf v}^-({\mathfrak g}(A_{r-1}^{(1)}))$ to the derived category 
which is used to construct $U_{\bf v}^-({\mathfrak g}(A_\infty))$. 
Let $\eta$ be anti-involutions on both quantum algebras which sends 
$f_i$ to $f_i$ respectively. 

Recall that ${\cal F}_{{\bf v}^{-1}}^{-\tilde\gamma}$ 
is a $U_{\bf v}({\mathfrak g}(A_\infty))$-module. We then have the 
following.

\begin{Theorem}[Varagnolo-Vasserot]

For each $x\in U_{\bf v}^-({\mathfrak g}(A_{r-1}^{(1)}))$, we set
\[
x\lambda=\sum_d \eta(\gamma_d(\eta(x)))
{\bf v}^{-\sum_{i<j,i\equiv j}d_ih_j}\lambda
\]
Then ${\cal F}_{{\bf v}^{-1}}^{-\tilde\gamma}$ becomes an 
$U_{\bf v}^-({\mathfrak g}(A_{r-1}^{(1)}))$-module. 
\end{Theorem}

\begin{Remark}
If we take $\tilde\gamma_i>\!\!>\tilde\gamma_{i+1}$, 
the canonical basis elements on these three Fock spaces 
coincide as long as the size of the Young diagrams 
indexing these canonical basis elements is not too large. 
\end{Remark}

\begin{Remark}
If we take $0\le-\tilde\gamma_1\le\cdots\le-\tilde\gamma_m<r$ in 
the above Fock space, we have 
Jimbo-Misra-Miwa-Okado higher level Fock space. This Fock 
space is the first example of higher level Fock spaces. 
\end{Remark}

By the above remark, we can use these Fock spaces to compute 
canonical basis elements on ${\cal F}_{\bf v}$ 
if we suitably care about the choice 
of $\tilde\gamma$. 

\begin{Theorem}[Uglov]
The Takemura-Uglov Fock space has a bar operation such that 
$\overline\emptyset=\emptyset$, 
$\overline{f_i\lambda}=f_i\overline\lambda$ and $\overline\lambda$ 
has the form $\lambda+$(higher terms) with respect to a dominance 
order. 
\end{Theorem}

The relation between the dominance order in the above theorem and 
the dominance order we use is well understood by using 
"abacus". As a conclusion, we can prove that 
$G(\lambda)=\lambda+\sum_{\mu\triangleright\lambda}c_{\lambda\mu}({\bf v})\mu$ 
as desired. 

We have explained that how crystal base theory on higher level 
Fock spaces fits in the modular representation theory of 
cyclotomic Hecke algebras. In particular, Kac $q$-dimension 
formula gives the generating function of the number of 
simple ${\cal H}_n$-modules. Even for type B Hecke algebras, 
it was new. 

\subsection{Future direction and Brou\'e's dream}
The original motivation of Brou\'e and Malle to introduce 
cyclotomic Hecke algebras is the study of modular representation 
theory of finite classical groups of Lie type over fields 
of non-describing characteristics. For example, Geck, Hiss and 
Malle's result towards classification of simple modules 
inspires many future problems. I may mention more in the lecture 
on demand. 

I would like to end these lectures with Brou\'e's famous dream. 
Let $B$ be a block of a group ring of a 
finite group $G$, and assume that it has an abelian 
defect group $D$. Let $b$ be the Brauer correspondent 
in the group ring of $DC_G(D)=C_G(D)\subset N_G(D)$. 
($(D,b)$ is called a maximal subpair or Brauer pair.) 
Then he conjectures that $D^b(B\!\!-\!\!mod)\simeq 
D^b(N_G(D,b)b\!-\!\!mod)$,i.e. 
$B$ and $N_G(D,b)b$ are derived equivalent
(Rickard equivalent). To be more precise on its base ring, let $(K,R,k)$ be 
a modular system. He conjectures the derived equivalence 
over $R$. 

Let $q$ be a power of a prime $p$, $G=G(q)$ be the general 
linear group 
$GL(n,q)$, and $k$ be an algebraically closed field of 
characteristic $l\ne p$, $(K,R,k)$ be a $l$-modular 
system.  Assume that $l>n$, and take 
$d$ such that $d|\Phi_d(q)$, 
$\Phi_d(q)| q^{n(n-1)/2}(q^n-1)\cdots(q-1)=|G(q)|$, where 
$\Phi_d({\bf q})$ is a cyclotomic polynomial. 
We take $B$ to be a unipotent 
block. In this case, unipotent blocks are 
paramerized by $d$-cuspidal pairs $(L(q),\lambda)$ up to 
conjugacy. Here $L(q)$ is a Levi subgroup, $\lambda$ is 
an irreducible cuspidal $KL(q)$-module. 
Further, $D$ is the $l$-part of the center of 
$L(q)$. ($L(q)$ is the centralizer of 
a "$\Phi_d$-torus" $S(q)$.) 
If we set $W(D,\lambda):=N_G(D,\lambda)/C_G(D)$, 
it is isomorphic to $G(d,1,a)$ for some $a$. 
$W(D,\lambda)$ is called {\bf cyclotomic Weyl group}. 
These are due to Brou\'e, Malle and Michel. 

In this setting, Brou\'e, Malle and Michel give an explicit 
conjecture on the bimodule complex which induces 
the Rickard equivalence between $B$ and $N_G(D,b)b$. 
It is given in terms of a variety which appeared in 
Deligne-Lusztig theory to trivialize 
a $L(q)$-bundle on a Deligne-Lusztig variety. 
Going down to the Deligne-Lusztig variety itself, 
it naturally conjectures the existence of a bimodule 
complex which induces derived equivalence between 
$B$ and a deformation ring 
of the group ring of the semi-direct of $S(q)_l$ 
with $W(D,\lambda)\simeq G(d,1,a)$. This conjecture 
is supported by the fact that they are isotypic 
in the sense of Brou\'e. 

It is expected that the deformation of $W(D,\lambda)$ is 
the cyclotomic Hecke algebra we have studied in these 
lectures. Hence, we expect that cyclotomic Hecke algebras 
with $m$ not restricted to $1$ or $2$ will have applications 
in this field. 
We remark that the Brou\'e conjecture is not 
restricted to $GL(n,q)$ only.

\end{document}